\documentclass{gtart}

\def\ifplaintex{\expandafter\ifx\csname documentclass\endcsname\relax}

\def\gtp{{\mathsurround=0pt\it $\cal G\mskip-2mu$eometry \&\ 
$\cal T\!\!$opology $\cal P\!$ublications}}  

\def\recd{{\small Received:\qua\receiveddate\ifx\reviseddate\relax
\else\qquad Revised:\qua\reviseddate\fi\par}} 

\def\lognumber#1{\def\thelognumber{#1}}
\def\volumenumber#1{\def\thevolumenumber{#1}}
\def\volumeyear#1{\def\thevolumeyear{#1}}
\def\papernumber#1{\def\thepapernumber{#1}}
\def\pagenumbers#1#2{\def\startpage{#1}\def\finishpage{#2}}
\def\published#1{\def\publishdate{#1}}

\def\received#1{\def\receiveddate{#1}}
\def\revised#1{\def\reviseddate{#1}}
\def\accepted#1{\def\accepteddate{#1}}

\long\def\asciiabstract#1{\long\def\theasciiabstract{#1}}

\let\\\par\let\thelognumber\relax\let\thevolumenumber\relax
\let\thepapernumber\relax\let\thevolumeyear\relax\let\startpage\relax
\let\finishpage\relax\let\publishdate\relax\let\receiveddate\relax
\let\reviseddate\relax\let\accepteddate\relax\let\theasciititle\relax
\let\theasciiauthors\relax
\let\theasciiabstract\relax

\let\theasciiemail\relax

\ifplaintex
\font\logobig=cmssbx10 scaled 3836
\font\logomed=cmssbx10 scaled 2557
\else
\font\logobig=cmssbx10 scaled 4200
\font\logomed=cmssbx10 scaled 2800
\fi

\long\def\makeagttitle{   
\count0=\startpage
\agt\hfill      
\hbox to 45truept{\vbox to 0pt{\vglue -13truept{\logomed A\kern -.37em{\logobig 
T}\kern -.38em G}\vss}\hss}
\break
{\small Volume \thevolumenumber\ (\thevolumeyear)
\startpage--\finishpage\nl
Published: \publishdate}

\vglue .25truein

{\parskip=0pt\leftskip 0pt plus
1fil\def\\{\par\smallskip}{\Large\bf\thetitle}\par\medskip} \vglue
0.05truein

{\parskip=0pt\leftskip 0pt plus 1fil\def\\{\par}{\sc\theauthors}
\par\medskip}%
 
\vglue 0.03truein 

{\small\leftskip 25truept\rightskip 25truept{\bf Abstract}\stdspace\theabstract

{\bf AMS Classification}\stdspace\theprimaryclass
\ifx\thesecondaryclass\relax\else; \thesecondaryclass\fi\par
{\bf Keywords}\stdspace \thekeywords\par}\vglue 7truept

}   

\ifplaintex
\font\phead=cmsl9 scaled 950
\font\pnum=cmbx10 scaled 913
\font\pfoot=cmsl9 scaled 950
\headline{\vbox to 0pt{\vskip -4.5mm\line{\small\phead\ifnum
\count0=\startpage ISSN 1472-2739 (on-line) 1472-2747 (printed)
\hfill {\pnum\folio}\else\ifodd\count0\def\\{ }%
\ifx\theshorttitle\relax\thetitle\else\theshorttitle\fi\hfill{\pnum\folio}
\else\def\\{ and }{\pnum\folio}\hfill\ifx\theshortauthors\relax\theauthors
\else\theshortauthors\fi\fi\fi}\vss}}
\footline{\vbox to 0pt{\vglue 0mm\line{\small\pfoot\ifnum\count0=\startpage
\copyright\ \gtp\hfill\else
\agt, Volume \thevolumenumber\ (\thevolumeyear)\hfill\fi}\vss}}
\else
\font=cmsl9 scaled 1050
\font\lnum=cmbx10 
\font=cmsl9 scaled 1050
\makeatletter
\def\@oddhead{{\small\lhead\ifnum\count0=\startpage ISSN 1472-2739 
(on-line) 1472-2747 (printed)\hfill {\lnum\number\count0}\else\ifodd\count0
\def\\{ }\ifx\theshorttitle\relax \thetitle \else\theshorttitle\fi\hfill
{\lnum\number\count0}\else\def\\{ and }{\lnum\number\count0}
\hfill\ifx\theshortauthors\relax 
\theauthors\else\theshortauthors\fi\fi\fi}}\def\@evenhead{\@oddhead}
\def\@oddfoot{\small\lfoot\ifnum\count0=\startpage\copyright\ \gtp\hfill\else
\agt, Volume \thevolumenumber\ (\thevolumeyear)\hfill\fi}
\def\@evenfoot{\@oddfoot}
\makeatother
\fi
\let\maketitlepage\makeagttitle
\let\makeshorttitle\maketitlepage
\let\maketitle\maketitlepage

\newwrite\gtoutfile
\long\gdef\makeheadfile{  
{\def\\{, }\def\s{ }
\immediate\openout\gtoutfile head.xxx
\immediate\write\gtoutfile{To: math@arxiv.org}
\immediate\write\gtoutfile{Subject: put OR rep NNNNN:ppppp}
\immediate\write\gtoutfile{--text follows this line--}
\immediate\write\gtoutfile{Proxy-for: \ifx\theasciiauthors\relax
\theauthors\else\theasciiauthors\fi\s<\ifx\theasciiemail\relax\theemail\else\theasciiemail\fi>}
\immediate\write\gtoutfile{\noexpand\\}
\immediate\write\gtoutfile{Authors: \ifx\theasciiauthors\relax
\theauthors\else\theasciiauthors\fi}
{\def\\{ }\immediate\write\gtoutfile{Title: \ifx\theasciititle\relax
\thetitle\else\theasciititle\fi}}
\immediate\write\gtoutfile{Subj-class: GT or SG, GR etc}
\immediate\write\gtoutfile{MSC-class: \theprimaryclass\ifx\thesecondaryclass\relax\else, \thesecondaryclass\fi}
\immediate\write\gtoutfile{Journal-ref: Algebr. Geom. Topol. \thevolumenumber\s
(\thevolumeyear) \startpage-\finishpage}
\immediate\write\gtoutfile{Comments: Published by Algebraic and
Geometric Topology at}
\immediate\write\gtoutfile{\s\s\s  http://www.maths.warwick.ac.uk/agt/AGTVol\thevolumenumber/agt-\thevolumenumber-\thepapernumber.abs.html}
\immediate\write\gtoutfile{\noexpand\\}
\immediate\write\gtoutfile{}
\ifx\theasciiabstract\relax
\immediate\write\gtoutfile{\theabstract}\else
\immediate\write\gtoutfile{\theasciiabstract}\fi
\immediate\write\gtoutfile{}
\immediate\write\gtoutfile{\noexpand\\}
\immediate\write\gtoutfile{}
\immediate\closeout\gtoutfile}}  

\def\maketitlepage{\makeagttitle\makeheadfile}
\let\makeshorttitle\maketitlepage
\let\maketitle\maketitlepage

\def\ifplaintex{\expandafter\ifx\csname documentclass\endcsname\relax}

\def\gtp{{\mathsurround=0pt\it $\cal G\mskip-2mu$eometry \&\ 
$\cal T\!\!$opology $\cal P\!$ublications}}  

\def\recd{{\small Received:\qua\receiveddate\ifx\reviseddate\relax
\else\qquad Revised:\qua\reviseddate\fi\par}} 

\def\lognumber#1{\def\thelognumber{#1}}
\def\volumenumber#1{\def\thevolumenumber{#1}}
\def\volumeyear#1{\def\thevolumeyear{#1}}
\def\papernumber#1{\def\thepapernumber{#1}}
\def\pagenumbers#1#2{\def\startpage{#1}\def\finishpage{#2}}
\def\published#1{\def\publishdate{#1}}

\def\received#1{\def\receiveddate{#1}}
\def\revised#1{\def\reviseddate{#1}}
\def\accepted#1{\def\accepteddate{#1}}

\long\def\asciiabstract#1{\long\def\theasciiabstract{#1}}

\let\\\par\let\thelognumber\relax\let\thevolumenumber\relax
\let\thepapernumber\relax\let\thevolumeyear\relax\let\startpage\relax
\let\finishpage\relax\let\publishdate\relax\let\receiveddate\relax
\let\reviseddate\relax\let\accepteddate\relax\let\theasciititle\relax
\let\theasciiauthors\relax
\let\theasciiabstract\relax

\let\theasciiemail\relax

\ifplaintex
\font\logobig=cmssbx10 scaled 3836
\font\logomed=cmssbx10 scaled 2557
\else
\font\logobig=cmssbx10 scaled 4200
\font\logomed=cmssbx10 scaled 2800
\fi

\long\def\makeagttitle{   
\count0=\startpage
\agt\hfill      
\hbox to 45truept{\vbox to 0pt{\vglue -13truept{\logomed A\kern -.37em{\logobig 
T}\kern -.38em G}\vss}\hss}
\break
{\small Volume \thevolumenumber\ (\thevolumeyear)
\startpage--\finishpage\nl
Published: \publishdate}

\vglue .25truein

{\parskip=0pt\leftskip 0pt plus
1fil\def\\{\par\smallskip}{\Large\bf\thetitle}\par\medskip} \vglue
0.05truein

{\parskip=0pt\leftskip 0pt plus 1fil\def\\{\par}{\sc\theauthors}
\par\medskip}%
 
\vglue 0.03truein 

{\small\leftskip 25truept\rightskip 25truept{\bf Abstract}\stdspace\theabstract

{\bf AMS Classification}\stdspace\theprimaryclass
\ifx\thesecondaryclass\relax\else; \thesecondaryclass\fi\par
{\bf Keywords}\stdspace \thekeywords\par}\vglue 7truept

}   

\ifplaintex
\font\phead=cmsl9 scaled 950
\font\pnum=cmbx10 scaled 913
\font\pfoot=cmsl9 scaled 950
\headline{\vbox to 0pt{\vskip -4.5mm\line{\small\phead\ifnum
\count0=\startpage ISSN 1472-2739 (on-line) 1472-2747 (printed)
\hfill {\pnum\folio}\else\ifodd\count0\def\\{ }%
\ifx\theshorttitle\relax\thetitle\else\theshorttitle\fi\hfill{\pnum\folio}
\else\def\\{ and }{\pnum\folio}\hfill\ifx\theshortauthors\relax\theauthors
\else\theshortauthors\fi\fi\fi}\vss}}
\footline{\vbox to 0pt{\vglue 0mm\line{\small\pfoot\ifnum\count0=\startpage
\copyright\ \gtp\hfill\else
\agt, Volume \thevolumenumber\ (\thevolumeyear)\hfill\fi}\vss}}
\else
\font=cmsl9 scaled 1050
\font\lnum=cmbx10 
\font=cmsl9 scaled 1050
\makeatletter
\def\@oddhead{{\small\lhead\ifnum\count0=\startpage ISSN 1472-2739 
(on-line) 1472-2747 (printed)\hfill {\lnum\number\count0}\else\ifodd\count0
\def\\{ }\ifx\theshorttitle\relax \thetitle \else\theshorttitle\fi\hfill
{\lnum\number\count0}\else\def\\{ and }{\lnum\number\count0}
\hfill\ifx\theshortauthors\relax 
\theauthors\else\theshortauthors\fi\fi\fi}}\def\@evenhead{\@oddhead}
\def\@oddfoot{\small\lfoot\ifnum\count0=\startpage\copyright\ \gtp\hfill\else
\agt, Volume \thevolumenumber\ (\thevolumeyear)\hfill\fi}
\def\@evenfoot{\@oddfoot}
\makeatother
\fi
\let\maketitlepage\makeagttitle
\let\makeshorttitle\maketitlepage
\let\maketitle\maketitlepage

\newwrite\gtoutfile
\long\gdef\makeheadfile{  
{\def\\{, }\def\s{ }
\immediate\openout\gtoutfile head.xxx
\immediate\write\gtoutfile{To: math@arxiv.org}
\immediate\write\gtoutfile{Subject: put OR rep NNNNN:ppppp}
\immediate\write\gtoutfile{--text follows this line--}
\immediate\write\gtoutfile{Proxy-for: \ifx\theasciiauthors\relax
\theauthors\else\theasciiauthors\fi\s<\ifx\theasciiemail\relax\theemail\else\theasciiemail\fi>}
\immediate\write\gtoutfile{\noexpand\\}
\immediate\write\gtoutfile{Authors: \ifx\theasciiauthors\relax
\theauthors\else\theasciiauthors\fi}
{\def\\{ }\immediate\write\gtoutfile{Title: \ifx\theasciititle\relax
\thetitle\else\theasciititle\fi}}
\immediate\write\gtoutfile{Subj-class: GT or SG, GR etc}
\immediate\write\gtoutfile{MSC-class: \theprimaryclass\ifx\thesecondaryclass\relax\else, \thesecondaryclass\fi}
\immediate\write\gtoutfile{Journal-ref: Algebr. Geom. Topol. \thevolumenumber\s
(\thevolumeyear) \startpage-\finishpage}
\immediate\write\gtoutfile{Comments: Published by Algebraic and
Geometric Topology at}
\immediate\write\gtoutfile{\s\s\s  http://www.maths.warwick.ac.uk/agt/AGTVol\thevolumenumber/agt-\thevolumenumber-\thepapernumber.abs.html}
\immediate\write\gtoutfile{\noexpand\\}
\immediate\write\gtoutfile{}
\ifx\theasciiabstract\relax
\immediate\write\gtoutfile{\theabstract}\else
\immediate\write\gtoutfile{\theasciiabstract}\fi
\immediate\write\gtoutfile{}
\immediate\write\gtoutfile{\noexpand\\}
\immediate\write\gtoutfile{}
\immediate\closeout\gtoutfile}}  

\def\maketitlepage{\makeagttitle\makeheadfile}
\let\makeshorttitle\maketitlepage
\let\maketitle\maketitlepage

\lognumber{14}
\volumenumber{1}
\volumeyear{2001}
\papernumber{14}
\pagenumbers{299}{310}
\received{9 October 2000}
\revised{11 May 2001}
\accepted{17 May 2001}
\published{23 May 2001}

\usepackage{amssymb}

\newtheorem{theorem}{Theorem}[section]
\newtheorem{lemma}[theorem]{Lemma}

\newtheorem{remark}[theorem]{Remark}

\newcommand{\tlk}{\mbox{\rm Tlk}}
\newcommand{\dlk}{\mbox{\rm Dlk}}

\newcommand{\cmapright}[2]{
\smash{\mathop{
\hbox to 1cm{\rightarrowfill}}\limits^{#1}_{#2}}}

\newcommand{\cmapleft}[2]{
\smash{\mathop{
\hbox to 1cm{\leftarrowfill}}\limits^{#1}_{#2}}}

\def\co{\colon\thinspace}

\begin{document}

\title{A Theorem of Sanderson on Link Bordisms\\in Dimension 4}

\authors{J. Scott Carter\\Seiichi Kamada\\Masahico Saito\\Shin Satoh}

\addresses{University of South Alabama, Mobile, AL 36688
\\Osaka City University, Osaka 558-8585, JAPAN
\\University of South Florida Tampa, FL 33620
\\RIMS, Kyoto University, Kyoto, 606-8502}                  
\email{carter@mathstat.usouthal.edu\\kamada@sci.osaka-cu.ac.jp\\saito@math.usf.edu\\satoh@kurims.kyoto-u.ac.jp}

\shorttitle{A theorem of Sanderson }
\shortauthors{Carter, Kamada, Saito and Satoh }

\begin{abstract}   

The groups of  link bordism
can be identified with homotopy groups
via the
Pontryagin--Thom construction.
B.J.~Sanderson
computed the bordism
group  of 3 component surface--links using the Hilton--Milnor Theorem,
and later gave a geometric interpretation of the groups
in terms of  intersections  of Seifert hypersurfaces
and their framings.
In this paper, we geometrically represent every element of  the bordism group
uniquely  by a certain standard form of a surface--link, a
generalization
of
a Hopf link. The standard forms give rise to  an inverse  of
Sanderson's geometrically defined invariant.

\end{abstract}

\asciiabstract{The groups of link bordism can be identified with
homotopy groups via the Pontryagin-Thom construction.  B.J. Sanderson
computed the bordism group of 3 component surface-links using the
Hilton-Milnor Theorem, and later gave a geometric interpretation of
the groups in terms of intersections of Seifert hypersurfaces and
their framings.  In this paper, we geometrically represent every
element of the bordism group uniquely by a certain standard form of a
surface-link, a generalization of a Hopf link. The standard forms
give rise to an inverse of Sanderson's geometrically defined
invariant.}


\maketitle  
%

\input epsf.tex

\def\figHopf2link{
\begin{figure}[ht!]
\begin{center}
\mbox{
\epsfxsize=3.6in
\epsfbox{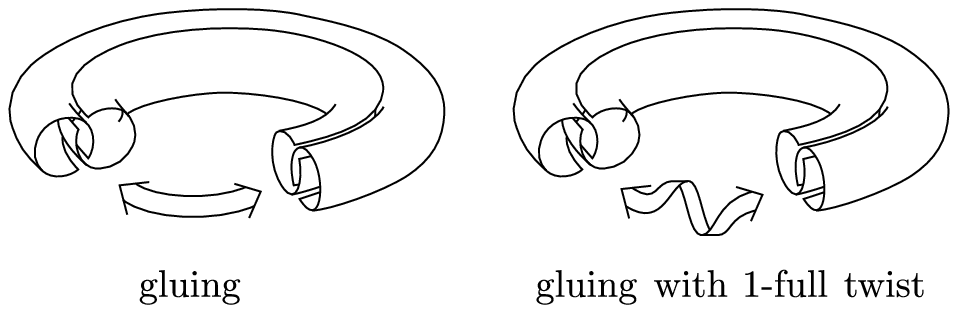}
}
\end{center}
\nocolon\caption{}
\end{figure}
}

\def\figPushOut{
\begin{figure}[ht!]
\begin{center}
\mbox{
\epsfxsize=4.7in
\epsfbox{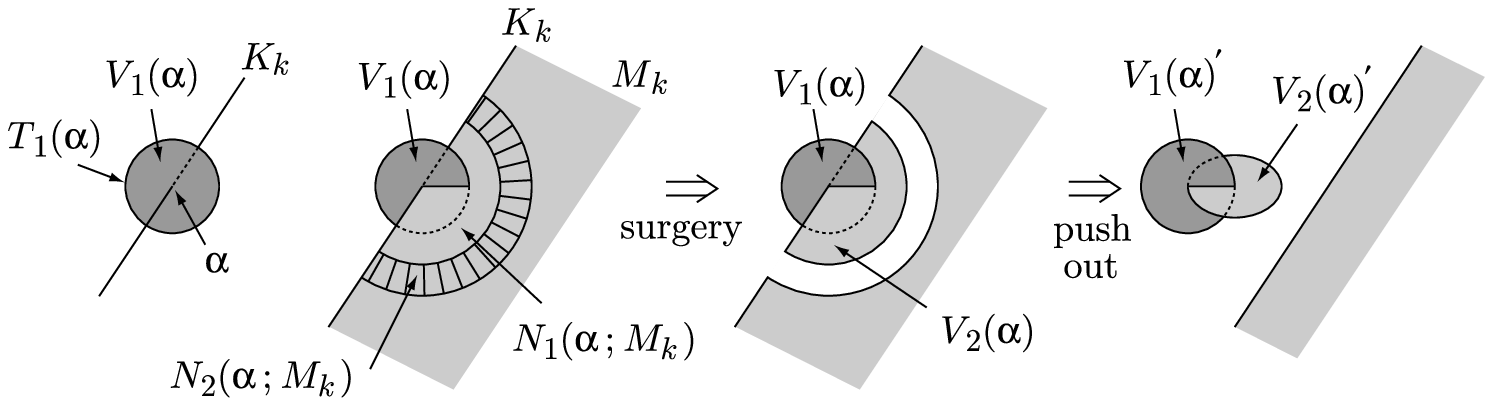}
}
\end{center}
\nocolon\caption{}
\end{figure}
}

\def\figFission{
\begin{figure}[htb]
\begin{center}
\mbox{
\epsfxsize=4in
\epsfbox{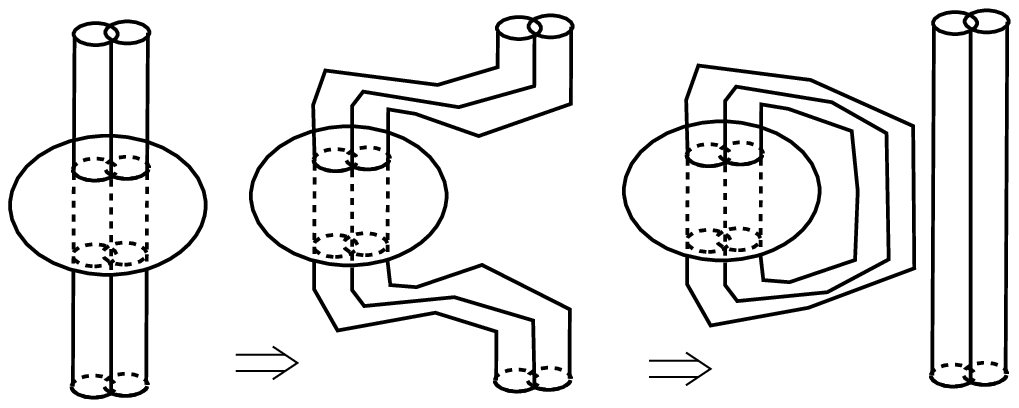}
}
\end{center}
\nocolon\caption{}
\end{figure}
}

\section{Introduction}

By  an {\it $n$-component surface--link}  we mean
a disjoint union  $F=K_1\cup \dots \cup K_n$ of
closed oriented surfaces $K_1, \dots, K_n$ embedded in ${\bf R}^4$
locally flatly.
Each $K_i$ is called
a {\it component}, which
may or may not be connected,
and possibly empty.
Two $n$-component surface--links
$F=K_1\cup \dots \cup K_n$ and $F'=K'_1\cup \dots \cup K'_n$
are {\it  bordant}
if there is a
compact oriented $3$-manifold $W$ properly embedded in  ${\bf R}^4
\times [0,1]$ such that $W$ has $n$ components $W_1, \dots, W_n$
with $\partial W_i= K_i \times \{0\} \cup (-K'_i)\times
\{1\}$.
Let $L_{4, n}$
    be the abelian group
of link bordism
classes of $n$-component surface--links.
The sum $[F] + [F']$ of $[F]$ and $[F']$
is defined by the class of
the split union of $F$ and $F'$.
The identity is represented by a
trivial $n$-component 2-link.
The inverse of $[F]$ is
represented by the mirror image of $F$ with the opposite
orientation.

\begin{sloppypar}
Brian Sanderson \cite{S:JLM}
identified the
bordism group $L_{4,n}$
of (oriented) $n$-component surface--links
with the homotopy group $\pi_4( \vee_{i=1}^{n-1}  S^2_i)$
using a generalized  Pon\-try\-a\-gin--Thom construction,
and computed the group
via the Hilton--Milnor Theorem.
Using the transversality theory developed in \cite{BRS},
Sanderson \cite{S:Bolyai} also
gave an explicit interpretation of the invariants
using Seifert hypersurfaces.
\end{sloppypar}

Specifically, let $F=K_1 \cup \cdots K_n$ be an $n$-component
surface--link. It is known
(cf.~\cite{NSato}) 
that each $K_i$ bounds a Seifert
hypersurface
$M_i$, i.e., an oriented compact $3$-manifold with
$\partial M_i = K_i$.
The {\it double linking number}, denoted by $\dlk(K_i, K_j)$, between
two components $K_i$ and $K_j$ is the framed intersection
$M_i \cdot K_j \in \pi_4(S^3) = {\bf Z}_2$.
The {\it triple  linking number}, denoted by $\tlk(K_i, K_j, K_k)$, among
three components $K_i$, $K_j$, and $K_k$, is the framed
intersection $M_i \cdot K_j \cdot M_k \in
\pi_4(S^4)
= {\bf Z}$.
Then Sanderson's geometrically defined invariant
$$ H \co L_{4, n} \rightarrow  A=(\underbrace{{\bf Z}\oplus\dots
\oplus{\bf Z}}_{\frac{n(n-1)(n-2)}{3}})
\oplus
(\underbrace{{\bf Z}_2 \oplus\dots
\oplus{\bf Z}_2}_{\frac{n(n-1)}{2}}) $$
is given by
$$ H([F])\!=\!( (\tlk(K_i, K_j, K_k),  \tlk(K_j, K_k, K_i)
)_{1 \leq i<j<k\leq n} ,
(\dlk(K_i, K_j) )_{1 \leq i<j \leq n} ) , $$
and it was shown \cite{S:Bolyai} that this gives an isomorphism.

The purpose of this paper is to give an inverse map
$$G \co A=(\underbrace{{\bf Z}\oplus\dots
\oplus{\bf Z}}_{\frac{n(n-1)(n-2)}{3}})
\oplus
(\underbrace{{\bf Z}_2 \oplus\dots
\oplus{\bf Z}_2}_{\frac{n(n-1)}{2}}) \rightarrow L_{4,n} $$
of $H$, by giving
an explicit set of
geometric representatives for a given value of Sanderson's invariant.
The representatives are
generalized Hopf links,
called
{\it Hopf $2$-links} (without or with {\it beads}),
which are
defined in
Sections~\ref{hopfsec} and \ref{beadsec}.

More specifically, we identify
$({\bf Z}\oplus\dots \oplus{\bf Z})\oplus
({\bf Z}_2 \oplus\dots \oplus{\bf Z}_2)$
with the abelian group
which
is abstractly generated by a certain family
${\cal F}$
of Hopf $2$-links
without or with beads and the homomorphism $G$ maps each generator
to its bordism class.
We prove that this homomorphism $G$ is
surjective.
(It is clear that  $H \circ G = \pm 1$; just use
the obvious Hopf solid link and/or normal 3-balls
as Seifert hypersurfaces when you evaluate $H$.)
Then $G$ is
an isomorphism which is an inverse of Sanderson's
homomorphism $H$
(up to sign).
The surjectivity is a consequence of the
following theorem.

\begin{theorem} \label{mainthm}
Any $n$-component surface--link $F$ is bordant to
a disjoint union of  Hopf $2$-links without or with beads.
More precisely, $[F] \equiv 0$ modulo
$\langle {\cal F} \rangle$,
where
$\langle {\cal F} \rangle$
is the subgroup of $L_{4,n}$ generated by
the classes of elements of ${\cal F}$.
\end{theorem}

Our proof constructs bordisms to unions of Hopf $2$-links geometrically,
and thus gives a self-contained geometric proof of
the Sanderson's classification theorem, without using
the
Hilton--Milnor theorem.

The paper is organized as follows.
Section~\ref{hopfsec} contains the definition of a
Hopf $2$-link and its characterization
in  the
bordism group.
Hopf $2$-links with beads and their roles in the bordism group
are
given in
Section~\ref{beadsec}, which also contains
the definition of the family ${\cal F}$ and
a proof of
    the Theorem~\ref{mainthm}.
We describe alternate
definitions of $\dlk$ and $\tlk$ in Section~\ref{linksec}.

\section{Hopf 2-Links}\label{hopfsec}

A {\it Hopf disk pair\/} is a
pair of disks $D_1$ and $D_2$ in a
3-ball $B^3$ such that
there is a homeomorphism from $B^3$ to a 3-ball
$$\{(x,y,z) ~ | ~ x^2 + y^2 + z^2 \leq 9 \}$$
which maps $D_1$ to
$$\{(x,y,z) ~ | ~ x^2 + y^2  \leq 1, z=0 \}$$
and $D_2$ to
$$\{(x,y,z) ~ | ~ (y-1)^2 + z^2  \leq 1, x=0 \}$$
homeomorphically.  The boundary of such a disk pair
is a Hopf link in $B^3$.

A pair of solid tori $V_1, V_2$ in ${\bf R}^4$ is
called a {\it Hopf solid link\/} if there is an embedding
$f \co B^3 \times S^1 \to {\bf R}^4$ such that
$f(D_i \times S^1) = V_i$ for $i=1,2$.
The boundary of a Hopf solid
link is called a {\it  Hopf $2$-link\/}, which is
a pair of embedded tori in  ${\bf R}^4$.
A simple loop $f({\rm (a~ point ~ of}~
B^3) \times S^1)$  is called a {\it core loop\/}
of the Hopf solid link and the Hopf $2$-link.
Any simple loop $\alpha$ in ${\bf R}^4$ is
ambient isotopic to a standard circle
in ${\bf R}^3 \subset {\bf R}^4$.
Since there are only two equivalence classes of
framings of $\alpha$ (or trivialization of
$N(\alpha) \cong B^3 \times S^1$),
any Hopf $2$-link is deformed by an ambient
isotopy  of  ${\bf R}^4$ so that the projection is one of
the illustrations depicted in
Figure~1.  If it is deformed into the
illustration on the
    left side, it is called a {\it standard Hopf $2$-link\/};
and if it is deformed into
the illustration on the
right side, it is called a {\it twisted Hopf $2$-link\/}.

We assume that a Hopf disk pair
is oriented so that the boundary is a positive Hopf link.
If the core loop is oriented, a Hopf solid link and a
Hopf $2$-link are assumed to be oriented by use of the orientation
of the Hopf disk pair and the orientation of the core loop.
(In this situation, we say that the Hopf $2$-link is oriented
coherently with respect to the orientation of the core loop.)

\begin{remark}\label{remark:hopf}
{\rm
Let $F=T_1 \cup T_2$ be a Hopf $2$-link whose core loop is
$\alpha$.  Then $F$ is ambient isotopic to
$-F = -T_1 \cup -T_2$ in $N(\alpha)$, where $-F$ means
$F$ with the opposite orientation.  Let $F'= -T_1 \cup T_2$
and $-F'= T_1 \cup -T_2$.  Then $F'$ and $-F'$ are ambient
isotopic in $N(\alpha)$.  When $\alpha$ is oriented,
one of $F$ and $F'$ is oriented coherently with respect to
$\alpha$, and the other is oriented coherently with respect to
$-\alpha$.
Since $\alpha$ and $-\alpha$ are ambient isotopic in ${\bf R}^4$,
$F$, $-F$, $F'$ and $-F'$ are ambient isotopic in ${\bf R}^4$.
}\end{remark}

\figHopf2link

\begin{lemma} \label{lem:standardHopf}
For a Hopf $2$-link $F= T_1 \cup T_2$, the following
conditions are mutually equivalent.
\begin{itemize}
\setlength{\itemsep}{3pt}
\item[{\rm (1)}]
$F$ is standard.
\item[{\rm (2)}]
$\dlk(T_1, T_2)=0$.
\item[{\rm (3)}]
$F$ is null-bordant.
\end{itemize}
\end{lemma}
\begin{proof}
Using a Hopf solid link, we see that
for the left side of Figure~1, $\dlk(T_1, T_2)=0$
and for the right,
$\dlk(T_1, T_2)=1$.  (This is also seen by
Remark~\ref{remark:dlk}.)
Thus (1) and (2) are equivalent.
Suppose (1).
Attach  $2$-handles to $T_1$ and $T_2$.
Then $T_1$ and $T_2$ change to $2$-spheres which
split by isotopy, and hence $F$ is null-bordant.
Thus (1) $\Rightarrow$ (3).
It is obvious that
(3) $\Rightarrow$ (2).
\end{proof}

Let $\alpha$ and $\alpha'$ be mutually disjoint
oriented simple loops in ${\bf R}^4$, and
let $F= T_1 \cup T_2$ and $F'= T_1' \cup T_2'$
be Hopf $2$-links whose
core loops are $\alpha$ and $\alpha'$
such that
$F$ and $F'$ are oriented coherently
with respect to the orientations
of $\alpha$ and $\alpha'$.
Let $\alpha''$ be an oriented loop obtained from
$\alpha \cup \alpha'$ by surgery along a band $B$
attached to $\alpha \cup \alpha'$.
Let $E$ be a $4$-manifold in ${\bf R}^4$ whose
interior contains $\alpha, \alpha'$ and $B$.

\begin{lemma}\label{lem:2compo}
In the above situation, there is
a Hopf $2$-link $F''= T_1'' \cup T_2''$
whose core loop is $\alpha''$ such that
$F''$ is bordant in $E$
to the $2$-component surface--link
$(T_1\cup T_1') \cup (T_2 \cup T_2')$.
\end{lemma}
\begin{proof}
Consider a tri-punctured sphere $\Sigma$ embedded
in $E \times [0,1]$
whose boundary is
$\partial \Sigma = \partial_0 \Sigma \cup
(-\partial_1 \Sigma)$ with
$\partial_0 \Sigma = (\alpha \cup \alpha') \times \{0\}$
and
$\partial_1 \Sigma = \alpha'' \times \{1\}$.
There exists an identification of
a regular neighborhood $N(\Sigma)$
in $E \times [0,1]$ with
$B^3 \times \Sigma$
such that
the Hopf $2$-links $F$ and $F'$ in $E \times \{0\}$
correspond to
$(\partial D_1 \cup \partial D_2) \times \partial_0 \Sigma$,
where $D_1 \cup D_2 \subset B^3$  is an oriented Hopf disk pair.
The desired $F''$ is obtained as
$(\partial D_1 \cup \partial D_2) \times \partial_1 \Sigma$.
\end{proof}

The transformation described in Lemma~\ref{lem:2compo}
is called {\it fusion} between
two Hopf $2$-links.
The inverse operation of fusion
is called {\it fission} of
a Hopf $2$-link.

\begin{lemma}\label{lem:2order}
Let $F= T_1 \cup T_2$ be a twisted Hopf $2$-link.
The order of $[F] \in L_{4,2}$ is two.
\end{lemma}
\begin{proof}
It is a consequence of Lemmas~\ref{lem:standardHopf} and
\ref{lem:2compo}.
\end{proof}

\section{Hopf $2$-Links with Beads and Proof of Theorem~\ref{mainthm} }
\label{beadsec}

Let $\alpha$ be a simple loop in ${\bf R}^4$ and let
$N(\alpha) \cong B^3 \times \alpha$ be a regular neighborhood.
We call a 3-disk $B^3 \times \{*\}$ ($* \in \alpha$) a
{\it meridian $3$-disk\/} of $\alpha$, and the boundary a
{\it meridian $2$-sphere\/} of $\alpha$.
    Let $D_1 \cup D_2$ be a Hopf disk pair in a $3$-disk $B^3$.
Let $f \co B^3 \times S^1 \to {\bf R}^4$ be an embedding, and let
$p_1, \dots, p_m$  be points of $S^1$.
We call the image $f(\partial D_1 \times S^1) \cup
f(\partial D_2 \times S^1) \cup
f( \partial B^3 \times \{p_1\}) \cup \dots \cup
f( \partial B^3 \times \{p_m\})$
a {\it Hopf $2$-link with beads\/}.  Each
meridian
$2$-sphere
$f( \partial B^3 \times \{p_i\})$ is called a {\it bead\/}.
We denote by  $S_{(i,j)}$  a twisted Hopf $2$-link
(as an $n$-component surface--link) at the
$i$\/th and the
$j$\/th component, and
by
$S_{(i,j,k)}$  a standard Hopf $2$-link
(as an $n$-component surface--link)
at the
$i$\/th and the
$j$\/th components  with a bead at the
$k$\/th component, respectively.
Let ${\cal F}$ denote a family of Hopf 2-links
$$
\{ S_{(i,j,k)} \,|\, i<j<k \} \cup
\{ S_{(i,k,j)} \,|\, i<j<k \} \cup
\{ S_{(i,j)} \,|\, i<j \}, $$
and $\langle{\cal F}\rangle$ the subgroup of
$L_{4,n}$ generated by the classes of elements of ${\cal F}$.

\begin{remark}\label{remark:hopfbead}
{\rm
Let $F= T_1 \cup T_2 \cup S$ be a Hopf $2$-link
$T_1 \cup T_2$
with a bead $S$ whose core is $\alpha$.
By Remark~\ref{remark:hopf},
$-T_1 \cup T_2 \cup S$ and $T_1 \cup -T_2 \cup S$
are ambient isotopic in $N(\alpha)$.  By an ambient
isotopy of ${\bf R}^4$ carrying $\alpha$ to $-\alpha$,
$-T_1 \cup T_2 \cup S$ is carried to
$T_1 \cup T_2 \cup -S$.  Therefore, any
Hopf $2$-link with a bead obtained from $F$ by
changing orientations of some components is
ambient isotopic to $F= T_1 \cup T_2 \cup S$
or $T_1 \cup T_2 \cup -S$.
}
\end{remark}

We denote by $S^-_{(i,j,k)}$ the Hopf 2-link $S_{(i,j,k)}$
such that the orientation of the bead is reversed.
It is clear that $[S^-_{(i,j,k)}] = -[S_{(i,j,k)}]$ in
$L_{4,n}$ (for example, use Lemmas~\ref{lem:standardHopf}
and \ref{lem:2compo}).
Thus, in the following proof,  we do not need to take care of an
orientation given to $S_{(i,j,k)}$.

\begin{proof}[Proof of Theorem~\ref{mainthm}]
(Step 1)
Let $M_1$ be a Seifert hypersurface for $K_1$
which intersects $K_2, \dots, K_n$ transversely.
For each $k$ ($k=2, \dots, n$), let $A_{1k}$
be the intersection $M_1 \cap K_k$, which is
the union of oriented simple loops in $K_k$ (or empty).
    Let $N(K_k)= D^2 \times K_k$ be a regular neighborhood
of $K_k$ in ${\bf R}^4$.
    For a component $\alpha$ of $A_{1k}$, let
$V_1(\alpha)$ be
a solid torus $D^2 \times \alpha$ $(\subset D^2 \times K_k =
N(K_k))$  in ${\bf R}^4$ and $T_1(\alpha)$ the boundary of
$V_1(\alpha)$.  They are oriented by use of the orientation of a
meridian disk
$D^2 \times \{*\}$ of $N(K_k)$ and the orientation of $\alpha$.
Then
$[\alpha] \in H_1(K_k)$ corresponds to
$[T_1(\alpha)] \in H_2(E(K_k))$ by the isomorphism
$$
H_1(K_k) \cong H^1(K_k) \cong
H^1(N(K_k)) \cong H_2(E(K_k))$$
obtained by the Poincar\'e and Alexander dualities,
where $E(K_k)$ is the exterior of $K_k$.
    Put
$V_1(A_{1k})= \bigcup_{\alpha \in A_{1k}} V_1(\alpha)$ and
$T_1(A_{1k})= \bigcup_{\alpha \in A_{1k}} T_1(\alpha)$.
    The surface $K_1$ is bordant to $\bigcup_{k=2}^n T_1(A_{1k})$
in ${\bf R}^4 \backslash (\bigcup_{k=2}^n K_k)$, for they cobound
a
3-manifold
${\rm Cl}(M_1 \backslash (\bigcup_{k=2}^n V_1(A_{1k}) ))$,
where ${\rm Cl}$ denotes the closure.
Thus, without loss of generality, we may assume that $K_1$ is
$\bigcup_{k=2}^n T_1(A_{1k})$.

Consider a Seifert hypersurface $M_k$ for $K_k$.
For a component $\alpha \in A_{1k}$,
let $V_2(\alpha) = N_1(\alpha; M_k)$ be a regular
neighborhood of $\alpha$ in $M_k$ such that
the union $V_1(\alpha) \cup V_2(\alpha)$ forms
a Hopf solid link in ${\bf R}^4$ with core $\alpha$,
see Figure~2 (the figure shows a section transverse to
$\alpha$).
      Let $N_2(\alpha; M_k)$ be a regular neighborhood of $\alpha$
in $M_k$ with $N_1(\alpha; M_k) \subset
{\rm int} N_2(\alpha; M_k)$ and let $C(\alpha)$ be a
3-manifold ${\rm Cl}(N_2(\alpha; M_k)\backslash N_1(\alpha; M_k))$.
By $C(\alpha)$, the $2$-component surface--link
$T_1(\alpha) \cup K_k$
is bordant to
$T_1(\alpha) \cup (\partial V_2(\alpha) \cup \partial
({\rm Cl}(M_k \backslash N_2(\alpha; M_k)))$.
This $2$-component surface--link is ambient isotopic
to
$\partial V_1(\alpha)' \cup (\partial V_2(\alpha)' \cup K_k)$,
where $V_1(\alpha)' \cup V_2(\alpha)'$ is a
Hopf solid link obtained from the Hopf solid link
$V(\alpha) \cup V_2(\alpha)$ by pushing out along
$N_1(\alpha; M_k)$ using $C(\alpha)$, see Figure~2.
      We denote by $\alpha'$ a loop obtained from $\alpha$
by pushing off along $M_k$ so that $\alpha'$ is
disjoint from $M_k$ and it is a
core of the Hopf solid link
$V_1(\alpha)' \cup V_2(\alpha)'$.

\figPushOut

Put
$V_1(A_{1k})'\!= \bigcup_{\alpha\in A_{1k} }\! V_1(\alpha)'$,
$V_2(A_{1k})'\!= \bigcup_{\alpha\in A_{1k} }\! V_2(\alpha)'$, and
$A_{1k}'\!= \bigcup_{\alpha \in A_{1k} }\! \alpha'$.
The $n$-component surface--link
$F= K_1 \cup \dots \cup K_n$ is bordant to
$F^{(1)}= K_1^{(1)} \cup \dots \cup K_n^{(1)}$ such that

$$\left\{
\begin{array}{l}
K_1^{(1)}= \bigcup_{k=2}^n \partial V_1(A_{1k})', \\ \\[-8pt]
K_j^{(1)}= \partial V_2(A_{1j})' \cup K_j \quad
{\rm for~} j=2, \dots n.
\end{array} \right.$$

(Step 2)
Let $M_2$ be the Seifert hypersurface for $K_2$ used in Step~1.
By the construction of $A_{12}'$, $M_2 \cap A_{12}' = \emptyset$.
For each $k$ with $3 \leq k \leq n$,
we may assume that $M_2$ intersects $K_k$ and
$A_{1k}'$ transversely.
Using $M_2$, we apply a similar argument to Step~1
to modify $K_2$ up to  bordism.
For a component $\alpha$ of $A_{2k} = M_2 \cap K_k$, let
$V_1(\alpha)$ be
a solid torus $D^2 \times \alpha$ $(\subset D^2 \times K_k =
N(K_k))$  in ${\bf R}^4$ and
$V_2(\alpha)$ a solid torus
$N_1(\alpha; M_k)$.
The intersection $M_2 \cap N(K_k)$
is $V_1(A_{2k}) = \bigcup_{\alpha \in A_{2k}} V_1(\alpha)$ and
the intersection   $M_2 \cap N(A_{1k}')$
is a union of some meridian 3-disks of $A_{1k}'$.
Let
$F^{(2)}$
be an $n$-component surface--link
$K_1^{(2)} \cup \dots \cup K_n^{(2)}$ such that
$$\left\{
\begin{array}{l}
K_1^{(2)}= \bigcup_{k=2}^n \partial V_1(A_{1k})', \\ \\[-8pt]
K_2^{(2)}= \partial V_2(A_{12})'
         \cup (\bigcup_{k=3}^n \partial V_1(A_{2k})'), \\ \\[-8pt]
K_j^{(2)}=
\partial V_2(A_{1j})' \cup \partial V_2(A_{2j})'
   \cup K_j \quad
{\rm for~} j=3, \dots n.
\end{array} \right.$$
By use of a 3-manifold which is $M_2$ removed the above
intersections, we see that
$F^{(1)}$ is bordant to
an $n$-component surface--link  $F^{(2)'}$  which is
the union of $F^{(2)}$
and some meridian $2$-spheres of $A_{1k}'$
with labels $2$
for
$k=3,\dots,n$.  The meridian $2$-spheres are the boundary
of meridian 3-disks that are
the intersection $M_2 \cap N(A_{1k}')$.
The surface--link $F^{(2)'}$ is bordant to
the union of $F^{(2)}$ and some standard Hopf $2$-links
$S_{(1,k,2)}$ or $S^-_{(1,k,2)}$
whose core loops are small trivial circles
in ${\bf R}^4$.  Figure~3 is a schematic picture of this
process (in projection in ${\bf R}^3$),
where an isotopic deformation
and fission of a Hopf $2$-link are applied.
Since $[S_{(1,k,2)}]$ belongs to
$\langle{\cal F}\rangle$,
we see that $F^{(2)'}$ is bordant to
$F^{(2)}$ modulo $\langle{\cal F}\rangle$.

\figFission

(Step 3)
Inductively, we see that $F^{(i-1)}$ is bordant to
$F^{(i)}= K_1^{(i)} \cup \dots \cup K_n^{(i)}$
modulo $\langle{\cal F}\rangle$
such that

$$
K_j^{(i)}=
\left\{
\begin{array}{ll}
(\bigcup_{k=1}^{j-1} \partial V_2(A_{kj})')
         \cup (\bigcup_{k=j+1}^n \partial V_1(A_{jk})')
\quad
&
{\rm for~} j {\rm ~with~} 1 \leq j \leq i,  \\ \\[-8pt]
(\bigcup_{k=1}^{i} \partial V_2(A_{kj})') \cup K_j
\quad
&
{\rm for~} j {\rm ~with~} i < j \leq n.
\end{array} \right.$$
Specifically, we see that
$F^{(i-1)}$ is bordant to an $n$-component
surface--link which is
a union of $F^{(i)}$
and some meridian 2-spheres of $A'_{k_1 k_2}$
such that $k_1 < i$ and $k_2 \neq i$.
These 2-spheres are the boundary of the meridian
3-disks that are the intersection
$M_i \cap N(A'_{k_1 k_2})$.  Since
$[S_{k_1, i, k_2}]$ belongs to
$\langle{\cal F}\rangle$,
$F^{(i-1)}$ is bordant to $F^{(i)}$ modulo
$\langle{\cal F}\rangle$.

Thus $F$ is bordant to
$F^{(n)}= K_1^{(n)} \cup \dots \cup K_n^{(n)}$
modulo $\langle{\cal F}\rangle$
such that
$$
\begin{array}{l}
K_j^{(n)}= (\bigcup_{k=1}^{j-1} \partial V_2(A_{kj})')
         \cup (\bigcup_{k=j+1}^n \partial V_1(A_{jk})')
\quad
{\rm for~} j =1, \dots,  n.
\end{array}
$$
It is a union of Hopf $2$-links and
the link bordism class is in $\langle{\cal F}\rangle$. 
\end{proof}

Recall that
Sanderson's homomorphism
$$H \co L_{4,n}\rightarrow A=(\underbrace{{\bf Z}\oplus\dots
\oplus{\bf Z}}_{\frac{n(n-1)(n-2)}{3}})
\oplus
(\underbrace{{\bf Z}_2 \oplus\dots
\oplus{\bf Z}_2}_{\frac{n(n-1)}{2}})$$
is defined by
$$ H([F])\!=( (\tlk(K_i, K_j, K_k),  \tlk(K_j, K_k, K_i)
)_{1 \leq i<j<k\leq n} ,
(\dlk(K_i, K_j) )_{1 \leq i<j \leq n})$$
for an $n$-component surface--link
$  F   =  K_1 \cup  \dots  \cup  K_n$.
Let $\{\   e_{ijk},\  e'_{ijk}\  | \   i<j<k \ \}\ \cup \ 
\{\ f_{ij}\ |\ i< j \ \}$ be a basis of $A$ such that
$e_{ijk}=(0,\dots,1,\dots,0)$
where $1$ corresponds to $\tlk(K_i,K_j,K_k)$,
$e_{ijk}'=(0,\dots,1,\dots,0)$
where $1$ corresponds to $\tlk(K_j,K_k,K_i)$,
and $f_{ij}=(0,\dots,1,\dots,0)$
where $1$ corresponds to $\dlk(K_i,K_j)$.
Give an orientation to each $S_{(i,j,k)}$
such that
$$
\left\{\begin{array}{ll}
\tlk(S_i, S_j, S_k) = +1, & \quad \tlk(S_k, S_j, S_i) = -1, \\[3pt]
\tlk(S_j, S_k, S_i) = 0, & \quad \tlk(S_i, S_k, S_j) = 0, \\[3pt]
\tlk(S_k, S_i, S_j) = -1, & \quad \tlk(S_j, S_i, S_k) = +1. \\[3pt]
\end{array}\right.
$$
We consider a homomorphism
$G \co A\rightarrow L_{4,n}$
with
$G(e_{ijk})=[S_{(i,j,k)}]$,
$G(e_{ijk}')= -[S_{(i,k,j)}]$, and
$G(f_{ij})=[S_{(i,j)}]$,
which is well-defined by Lemma~\ref{lem:2order}.
Then we have $H \circ G = {\rm id}$.
Theorem~\ref{mainthm} says that
$G$ is surjective.

\begin{remark}\label{remark:tlkeqs}
{\rm
In general, for any $n$-component surface--link $F$, it
is known (cf.~\cite{CJKLS})
that for any distinct $i,j,k$,
$$
\begin{array}{c}
\dlk(K_i, K_j)=\dlk(K_j, K_i),  \\[3pt]
\tlk(K_i, K_j, K_k)  =  -\tlk(K_k, K_j, K_i),   \\[3pt]
\tlk(K_i, K_j, K_k) + \tlk(K_j, K_k, K_i)+ \tlk(K_k, K_i, K_j)=0.
\end{array}
$$
{}From these formulas, it follows
immediately
from the condition
$H([F])=H([F'])$
that
$\dlk(K_i, K_j)\!=\!\dlk(K'_i, K'_j)$
and $\tlk(K_i, K_j, K_j)\!=\!\tlk(K'_i, K'_j, K'_k)$ for
any distinct $i,j,k$.
}\end{remark}

\section{Remarks on
Linking Numbers} \label{linksec}

In this section we comment on different definitions and
aspects of the generalized linking numbers, $\dlk$ and $\tlk$.
To mention an analogy to projectional definition
    of the classical linking number
(cf.~\cite{Rolfsen}),  
we start with a review of projections
of surface--links.

Let $F=K_1 \cup \cdots \cup K_n$ be an $n$-component
surface-- and let $F^*= p(F)$ be a surface diagram of $F$
with respect to a projection $p \co {\bf R}^4\rightarrow{\bf R}^3$.
For details of  the definition of a surface diagram, see
\cite{CS:Reidemeister, CS:book, Rose} 
for example.
The singularity set of $F^*$
consists of double points and
isolated branch/triple points.  The singularity set is
a union of immersed circles and arcs
in ${\bf R}^3$, which is called
the set of
{\it double curves}.
Two sheets intersect along a double curve,
which are called {\it upper} and {\it lower}
with respect to the projection direction.
A double curve is {\it of type $(i,j)$}
if the upper sheet comes from $K_i$ and
the lower comes from $K_j$.
We denote by $D_{ij}$ the union of
double curves of type $(i,j)$.
    For distinct $i$ and $j$, $D_{ij}$ is the union of
immersed circles.
(If $D_{ij}$ contains an immersed arc, its end-points
are branch points.  So the upper sheet with label $i$
and the lower
sheet with label $j$ along the arc come from the same
component of $F$.  This contradicts $i\neq j$.)
Let $D_{ij}^+$ be the union of
immersed circles in ${\bf R}^3$
obtained from $D_{ij}$ by shifting it
in a diagonal direction that is in
the positive normal direction of
the upper sheet and also
in the positive normal direction of the lower
sheet of
$F^*$  along $D_{ij}$
so that $D_{ij}$ and $D_{ij}^+$ are disjoint.
    Let  $\widetilde{D}_{ij}$ be a link (i.e.,
embedded circles)  in
${\bf R}^3$ which is obtained from $D_{ij}$
by a slight perturbation by a homotopy, and
let  $\widetilde{D}_{ij}^+$ be a link in
${\bf R}^3$ which is obtained from $D_{ij}^+$
similarly.
Give  $\widetilde{D}_{ij}$ an orientation
and $\widetilde{D}_{ij}^+$ the orientation
which is parallel to
that
of
$\widetilde{D}_{ij}$.
The linking number between $\widetilde{D}_{ij}$
and $\widetilde{D}_{ij}^+$ does not depend on
the perturbations and
the orientation of $\widetilde{D}_{ij}$, which we call
the {\it linking number} between $D_{ij}$
and $D_{ij}^+$.
Then
we have

\begin{remark}\label{remark:dlk}
{\rm
The  double linking number $\dlk(K_i, K_j)$ is equal to
a value in ${\bf Z}_2=\{0,1\}$ that is
the linking number between $D_{ij}$
and $D_{ij}^+$ modulo $2$.
}\end{remark}

    At a triple point in the projection, three
sheets intersect that have distinct relative heights with respect to the
projection direction, and we call them
{\it top}, {\it middle}, and {\it bottom} sheets,
accordingly.
If the orientation normals to the top, middle, bottom sheets at
a triple point
$\tau$ matches with this order the fixed orientation of ${\bf R}^3$,
then the sign of $\tau$ is positive and $\varepsilon(\tau)=1$.
Otherwise the sign is negative and $\varepsilon(\tau)=-1$.
    (See  \cite{CJKLS, CS:book}.)
A triple point is {\it of type $(i,j,k)$}
if the top sheet comes from $K_i$,
the middle comes from $K_j$, and
the bottom comes from $K_k$.
The following projectional interpretation of
triple linking numbers was extensively used in
\cite{CJKLS} for invariants defined from quandles.

\begin{remark} {\rm
The triple linking number $\tlk(K_i,K_j,K_k)$
is (up to sign) the sum of the signs
of all the triple points of type
$(i,j,k)$.}
\end{remark}

Let $f \co F_1 \cup F_2 \cup F_3 \rightarrow {\bf R}^4$ denote an
embedding of the disjoint union of oriented surfaces $F_i$
representing $F=K_1\cup K_2\cup K_3$.
Define a map
$L \co F_1 \times F_2 \times F_3 \rightarrow S^3 \times S^3$ by
$$L(x_1,x_2,x_3) = \left( \frac {f(x_1)-f(x_2)}{||f(x_1)-f(x_2)||},
\frac {f(x_2)-f(x_3)}{||f(x_2)-f(x_3)||} \right)$$
for $x_1\in F_1$, $x_2\in F_2$ and $x_3 \in F_3$.
In \cite{Koschorke0} it is observed that
the degree of $L$ is (up to sign)
the triple linking number $\tlk(K_1,K_2,K_3)$.

For further related topics, refer to
\cite{BT, Tim1, TO, FennRolf, Kirk, KK, Koschorke,  MasseyRolf, Ruberman}.

\medskip

\noindent
{\large\bf Acknowledgments} \
The authors would like to thank Akio Kawauchi,
    Brian Sanderson, and the referee
    for helpful suggestions.
JSC is being supported by NSF grant DMS-9988107.
SK and SS are being supported by Fellowships
from the Japan Society for the Promotion of Science.
MS is being supported by NSF grant DMS-9988101.

\Addresses\recd


\begin{thebibliography}{99}


\bibitem{BT} Bartels, A. and Teichner, P.,
{\it All two dimensional links are null homotopic,}
Geometry and Topology, Vol. 3 (1999),  235--252.

\bibitem{BRS}
    Buoncristiano, S., Rourke, C.P., Sanderson, B.J.,
{\it  A geometric approach to
homology theory,}
    London Mathematical Society Lecture Note Series, No. 18. Cambridge 
University
Press, Cambridge-New York-Melbourne, 1976.




\bibitem{CJKLS}
Carter,  J.S., Jelsovsky, D.,  Kamada, S.,
Langford, L., and Saito, M.,
{\it Quandle cohomology and state-sum invariants
of knotted curves and surfaces,}
preprint at
http://xxx.lanl.gov/abs/math.GT/9903135 .


\bibitem{CS:Reidemeister}
Carter, J.S. and  Saito, M.,
{\it Reidemeister moves for surface isotopies and
their interpretations as moves to movies,}
J. Knot Theory Ramifications 2 (1993), 251--284.

\bibitem{CS:book}  Carter, J.S. and  Saito, M.,
{\it Knotted surfaces and their diagrams,}
 American Mathematical Society, Mathematical Surveys and Monograph Series, Vol 55, (Providence 1998).


\bibitem{Tim1}
Cochran, T.D.,
{\it  On an invariant of link cobordism in dimension four,}
Topology Appl. 18 (1984), no. 2-3, 97--108.

\bibitem{TO}
Cochran, T.D. and  Orr, K.E.,
{\it  Not all links are concordant to boundary
links, }  Ann. of Math. (2) 138 (1993), no. 3, 519--554.

\bibitem{FennRolf}
Fenn, R.  and  Rolfsen, D.,
{\it Spheres may link homotopically in
$4$-space,}
    J. London Math. Soc. (2) 34 (1986), no. 1, 177--184.


\bibitem{Kirk} Kirk, P.A.,
{\it Link maps in the four sphere,}
Differential topology (Siegen, 1987), 31--43,
Lecture Notes in Math., 1350,
Springer, Berlin-New York, 1988.

\bibitem{KK}
   Kirk, P.A. and  Koschorke, U.,
{\it Generalized Seifert surfaces and linking numbers,}
Topology Appl. 42 (1991), no. 3, 247--262.

\bibitem{Koschorke}
   Koschorke, U.,
{\it  Homotopy, concordance and bordism of link maps,}
Global analysis in modern mathematics (Orono, ME, 1991; Waltham, MA, 1992),
283--299,
Publish or Perish, Houston, TX, 1993.

\bibitem{Koschorke0}
   Koschorke, U.,
{\it A generalization of Milnor's $\mu$-invariants
    to higher-dimensional link maps,}
Topology 36, 2 (1997), 301--324.



\bibitem{MasseyRolf}  Massey, W.S.  and  Rolfsen, D.,
{\it Homotopy classification of higher-dimensional links,}
Indiana Univ. Math. J. 34 (1985), no. 2, 375--391.

\bibitem{Rolfsen} Rolfsen, D.,
{\it Knots and links,}
Publish or Perish, Inc., 1976.

\bibitem{Rose} Roseman, D.,
{\it Reidemeister-type moves for surfaces in four dimensional space, }
in Banach Center Publications 42 (1998) Knot theory, 347--380.

\bibitem{Ruberman}
Ruberman, D.,
{\it  Concordance of links in $S\sp 4$,}
    Four-manifold
theory (Durham, N.H., 1982), 481--483, Contemp. Math., 35, Amer. Math.
Soc., Providence, R.I.,
1984.

\bibitem{S:JLM}
Sanderson, B. J.,
{\it  Bordism of links in codimension $2$,}
    J. London Math. Soc.
(2) 35 (1987), no. 2, 367--376.




\bibitem{S:Bolyai}
  Sanderson, B. J.,
{\it  Triple links in codimension $2$,}
  Topology. Theory and
applications, II (P\'{e}cs, 1989), 457--471, Colloq. Math. Soc. J\'{a}nos 
Bolyai, 55, North-Holland,
Amsterdam, 1993.




\bibitem{NSato} Sato, N.,
{\it Cobordisms of semiboundary links.}
Topology Appl. 18
    (1984), no. 2-3, 225--234.

\end{thebibliography}
\end{document}